\documentclass{article}

\title{Scaling of Percolation on Infinite Planar Maps, I}

\author{Omer Angel}

\date{November, 2004}

\usepackage{amsmath,amsthm,amssymb,amsfonts}
\usepackage{graphics,subfigure,pst-all}
\usepackage{hyperref}

\newtheorem{thm}{Theorem}[section]

\newtheorem{coro}[thm]{Corollary}
\newtheorem{prop}[thm]{Proposition}
\theoremstyle{definition}

\newcommand{\thmref}[1]{Theorem~\ref{T:#1}}
\newcommand{\corref}[1]{Corollary~\ref{C:#1}}
\newcommand{\figref}[1]{Figure~\ref{fig:#1}}

\newcommand{\Z}{\ensuremath{\mathbb{Z}}}

\newcommand{\R}{\ensuremath{\mathbb{R}}}
\newcommand{\E}{\ensuremath{\mathbb{E}}}
\renewcommand{\P}{\ensuremath{\mathbb{P}}}
\newcommand{\tP}{\ensuremath{\widetilde \P}}
\newcommand{\tA}{\ensuremath{\widetilde A}}
\newcommand{\tQ}{\ensuremath{\widetilde Q}}

\begin{document}

\maketitle
\begin{abstract}
  We consider several aspects of the scaling limit of percolation on random
  planar triangulations, both finite and infinite. The equivalents for
  random maps of Cardy's formula for the limit under scaling of various
  crossing probabilities are given. The limit probabilities are expressed
  in terms of simple events regarding Airy-L\'evy processes. Some explicit
  formulas for limit probabilities follow from this relation by applying
  known results on stable processes. Conversely, natural symmetries of the
  random maps imply identities concerning the Airy-L\'evy processes.
\end{abstract}

\section{Introduction}
%%%%%%%%%%%%%%%%%%%%%%

Random planar maps have long been of interest in both the mathematical
\cite{Tutte1, Tutte2, BRW, RiWo} and physical \cite{ADJ,AmWa} literature.
The commonly used distributions on planar maps are the uniform distribution
on maps of a given size, and Boltzmann distributions, where the weight of a
map is exponential in its size. In \cite{UIPT1}, another natural measure
was introduced --- the uniform infinite planar triangulation (UIPT). The
UIPT is a probability measure supported on infinite rooted planar
triangulations. The UIPT is uniquely characterized, in that it satisfies a
Gibbsian local uniformity property, roughly stating that disjoint regions
in the map are independent and each has a uniform distribution.

While the UIPT construction extends to many classes of planar maps, this
paper will concentrate on triangulations. The natural connection between
percolation and random triangulations is akin to the relation to the
Euclidean triangular lattice. This connection has been partially exploited in
\cite{UIPT2}, and is fundamental to this paper as well.

There are several variations on the definition of a triangulation, with the
most common difference being constrants on the connectivity. Within the
class of all triangulations, there is a smaller class of 2-connected
triangulations where self-loops are forbidden. Smaller still classes
contains 3-connected which have no multiple edges, and 4-connected
triangulations, where there is no triangle in the underlying graph which is
not a face of the map. Other sub-classes are based on restrictions on
vertex degrees, etc. In this paper we consider the case of 2-connected
planar triangulations since this class of maps is most amenable to the
tools introduced in \cite{UIPT2}. However, it should be noted that for most
questions concerning percolation on the triangulations, triangulations of
the different connectivity classes are in fact equivalent. Indeed, for
combinatorial reasons, the percolation restricted to the 4-connected core
(see \cite{BFSS}) of a triangulation contains most of the relevent
information.

\medskip

One of the key problems concerning random planar maps, which has been
partially solved in \cite{ChSc,MaMo} is to understand the scaling limit of
uniform planar maps. It is believed that there is some measure on random
metric spaces, which is a scaling limit of the uniform sample of planar
maps. Furthermore, this measure is canonical in that it does
not depend (up to some minimal restrictions) on the class of planar maps
used at the discrete level (as the Brownian motion is the scaling limit of
random walks with weak restrictions.)

Here we deal not with scaling of the map itself, but rather with the way
critical percolation on the map is scaled. In particular, we consider the
random map equivalent of Cardy's formula for the Euclidean plane. Start
with a region in the plane with two disjoint segments along its boundary,
and consider the probability that critical percolation in the domain
contains a crossing between the two segments. In the Euclidean plane this
probability is believed to have a scaling limit (as the lattice mesh size
tends to 0) which is a conformally invariant function of the domain.
Smirnov \cite{Smir} has shown this to be the case for site percolation on
the triangular lattice.

In the random map setting, the situation is somewhat simpler. The boundary
of a domain contains very little information, namely its size. No geometric
information is stored in the boundary. Equally, when the boundary has
several segments marked out on it, the only information relevant to the map
inside is the lengths of the segments and the distances between them along
the boundary. This is part of the reason that percolation on random maps is
easier to fully understand than in the Euclidean plane.

Given the (not fully understood) connection between random planar maps and
Euclidean lattices \cite{dupl1,dupl2}, a connection that is embodied by the
KPZ relation \cite{KPZ}, it is of interest to better understand how
critical models behave on random planar maps. In particular, it is of
interest to calculate critical exponents in the random map setting.

\medskip

In this paper we consider the crossing probabilities for critical
percolation on the UIPT. Planar maps come in several flavors: Finite and
infinite, and with a finite boundary, no boundary or an infinite one. These
are described below in detail. Of particular interest is the half plane
UIPT. As the name suggests this a random triangulation in the half plane,
i.e.\ having a boundary that is an infinite line. The half plane UIPT has
the remarkable property, reminiscent of SLE curves, that when a piece of it
is removed, the remaining triangulation has the same law as the whole
triangulation. 

The main results of this paper express the crossing probabilities in a
number of scenarios, in terms of hitting probabilities of a random walk
with i.i.d.\ steps with distribution
\[
  X_i = \begin{cases}
     1  & \text{with prob. $2/3$,} \\
    -k  & \text{with prob. $2\frac{(2k-2)!}{4^k(k-1)!(k+1)!}$.} \\
  \end{cases}
\]
(We later use a lazy version of this random walk which is more natural.
This has no effect on hitting probabilities.) See Theorems~\ref{T:2segs}
and \ref{T:3segs} for precise statements.

The simplest manifestation of the relation between crossing probabilities
and the random walk is given in \thmref{2segs}, and states that a certain
crossing probability $Q_{a,b}$, is given by
\[
  Q_{a,b} = \P_{a}(S \textrm{ hits $\Z^-$ in $(-\infty,-b]$}),
\]
where $a,b$ describe the boundary conditions, for the crossing problem, and
$S$ is a sample of the above random walk started at $a$.

Crossing probabilities in more complex scenarios can similarly be expressed
in terms of more complicated hitting probabilities of the same random walk,
or multiple independent copies of the random walk. In this paper we
consider several variations of crossing probabilities in the half plane and
in Boltzmann triangulations. In a subsequent paper we consider more general
events on percolation, on these triangulations as well as in the UIPT in a
disc.

A second theme of the paper is to consider the limit under scaling of the
crossing probabilities. Since the random walks have a scaling limit given
by a stable process with increment having the Airy distribution, the
crossing probabilities have a scaling limit given by hitting distributions
of such processes, see Theorems~\ref{T:sclim2}, \ref{T:sclim2a}

Finally, natural symmetries of the triangulations are expressed in terms of
several identities on the Airy-stable processes. The simplest of these
identities (\corref{symmetry}) states that if $\P_a$ is the probability
measure of an Airy-stable process $Y_t$ started at $a$, and hitting $\R^-$
at time $T_-$, then
\[
  \P_a\left(|Y_{T_-}| > b\right) = \P_b\left(|Y_{T_-}| < a\right).
\]
This identity can also be derived directly from known properties of the
process, but identities arising from more involved percolation problems are
more difficult to derive directly, see Corollaries~\ref{C:id2},
\ref{C:id3}, \ref{C:id4}.

\medskip

The next section reviews the UIPT, and introduces the half plane UIPT, and
recalls some properties of the Airy-Stable process.
Section~\ref{sec:half_plane} contains calculations of crossing
probabilities in the half plane UIPT, as well as their scaling limits.
Section~\ref{sec:boltzmann} uses the fundamental relation between infinite
triangulations and the Boltzmann distribution on finite triangulations to
calculate crossing probabilities on the latter.

\section{Background}
%%%%%%%%%%%%%%%%%%%%

\subsection{UIPT in the plane, half plane and disc}

When considering a random planar map, one starts with some class of planar
maps that is typically defined by some combinatorial constraints (such as
having faces or vertices of certain degrees, connectivity, etc.) However,
beyond the choice of class there are several possible large scale
topological types of planar maps, each with its associated uniform
distribution. The first and most studied is the map on the sphere, where
there is no external face. Since crossing probabilities require a boundary,
we will not be concerned with random maps on the sphere.

An important distribution on finite planar maps that we will study is the
{\em free} or {\em Boltzmann} distribution on maps in a disc. These
distributions have played a role in \cite{UIPT1}. We use the following
notations (all of which can be extended to other classes of planar maps):

\begin{description}
\item[$\phi_{n,m}$] The number of triangulations of the $m$-gon of size $n$
  (with $n$ internal vertices). These have asymptotics $\phi_{n,m} \approx
  C_m n^{-5/2} \alpha^n$ as $n\to\infty$.
\item[$\alpha$] The exponential growth rate of $\phi_{n,m}$, which happens
  to be 27/2.
\item[$C_m$] The boundary size component in the asymptotics of
  $\phi_{n,m}$. These are given by
  \[
    C_m = \frac{4(2m-3)!}{3^{7/2}\sqrt{\pi}(m-2)!^2} (9/4)^m.
  \]
\item[$Z_m$] The partition function for triangulations of an $m$-gon:
  \[
  Z_m = \sum \phi_{n,m} \alpha^{-n} = \frac{4(2m-4)!}{9m!(m-2)!} (9/4)^m.
  \]
\item[$\mu_m$] The free (Boltzmann) measure on triangulations of an
  $m$-gon, assigning to each triangulation $T$ a probability $\mu_m(T) =
  Z_m^{-1} \alpha^{-|T|}$. The free map is a sample from $\mu_m$.
\end{description}

Apart from measures on spheric maps and the free map, there are several
natural measures supported on infinite maps, each having a different global
geometry, and these will all be referred to as UIPT's. The distinction is a
discrete variation of the conformal distinction between the full complex
plane and domains with non-trivial complement. In the discrete setting,
there is also a distinction between domains with finite boundary and
domains with infinite boundary. These geometries are represented by the
full plane, the half plane and the disc. The full plane UIPT is the one
defined in \cite{UIPT1}. As an aside, note that there is no infinite
spheric map, since it must have an accumulation point somewhere on the
sphere. However, when taking a scaling limit of the triangulations, it is
possible to get a topological sphere.

The UIPT in a disc or polygon (since the number of vertices on the boundary
is fixed) is referred to --- though not explicitly defined --- in
\cite{UIPT2}. The UIPT of a polygon can be defined by starting with a
sample $T$ of the UIPT, conditioning the ball of radius 1 around the root
to be a wheel graph with $m$ spokes, (a cycle with edges to the root), and
then removing that wheel from $T$ and marking a new root on the resulting
face. When the map is embedded in the plane, inversion may be applied to
have a triangulation of the interior of a disc rather then exterior.

It follows as a corollary of the Locality Theorem of \cite{UIPT1}, that
replacing the wheel in the previous paragraph by any other simply connected
neighborhood of the root with boundary size $m$ results in the same
distribution. The same is true if one starts with a UIPT of an $m'$-gon and
conditioning on some structure in an annulus with boundaries of length $m'$
and $m$. Thus the resulting distribution is canonical and does not depend
on the choice of neighborhood.

The three UIPT measures shall be denoted as follows:

\begin{description}
\item[$\nu$] The full-plane UIPT of \cite{UIPT1}.
\item[$\nu_m$] The UIPT in an $m$-gon, which is also the limit as
  $n\to\infty$ of $\mu_m(\cdot \,\vert\,|T|>n)$, i.e.\ the limit of $\mu_m$
  conditioned on having a large volume.
\item[$\nu_\infty$] The UIPT in the half plane, formally defined below.
\end{description}

The UIPT in the half plane is of special interest, since the exploration
process described below takes a particularly simple form there. The half
plane UIPT may be defined by taking a limit of the UIPT measures in discs
with boundary sizes tending to infinity. This is a weak limit with respect
to the topology described in \cite{UIPT1}.

\begin{thm}\label{T:hp_limit}
  With the above notations, The measures $\nu_m$ have a weak limit
  $\nu_\infty$ as $m\to\infty$. The limit is a probability measure
  supported on infinite triangulations of the half plane. The weak limit as
  $m\to\infty$ of $\mu_m$ also exists and equals $\nu_\infty$ as well.
\end{thm}

The UIPT of the half plane may also be constructed by ``conditioning'' the
full plane or disc UIPT to include some infinite slit. This condition has
probability 0, so such a definition would also be stated as a weak limit of
measures, which turn out to the same sequence of approximations as in
\thmref{hp_limit}.

The proof is based on the peeling procedure, which is the main technique
developed in \cite{UIPT2}. The procedure in the disc and full plane UIPT is
described below for completeness, though we recast it here as an
exploration process, where initially the random map is unknown and at each 
step some finite portion of the map is revealed. Formally this is defined
as a filtration of the UIPT probability space. 

The exploration process of a triangulation (or general map) proceeds as
follows: A boundary edge is chosen (randomly, or in a deterministic manner
depending on the past), and the triangle (or face) containing the edge is
revealed. In the case of a triangulation there are two possibilities.
Either the third vertex is some other boundary vertex, or it is an internal
vertex. These events have probabilities depending on the measure used. In
the case of more general maps, there might of course be more possibilities.

For each of the measures mentioned above (i.e.\ $\mu_m$, $\nu$, $\nu_m$, or
$\nu_\infty$) the following holds: If the third vertex was an internal
vertex, then the remaining unexplored part of the map has the same
distribution (with $m+1$ for $m$ when the boundary is finite). If the third
vertex was a boundary vertex at distance $k$ from the edge used, then the
unexplored part of the map contains two components, one with boundary size
$k+1$, and the other with boundary size $m-k$ or $\infty$ as the case may
be. These two components are independent, one of them contains a free
triangulation of a disc, and the other contains the same measure as the
original map that is explored (again, with different boundary size if
applicable). In the case of $\mu_m$, this means both components have the
free distribution.

In the case of exploration of an infinite map, the exploration process
never terminates. Of course, some events may be in the sigma algebra
generated by some finite time, and to determine their occurrence there is no
need to continue the exploration process. In the infinite cases, when a
component is generated that contains a free triangulation, one can either
reveal that triangulation as part of the same step of the exploration
process, or leave it unknown. For our purposes it is more convenient to
leave it as unknown. If such a component needs to be further explored, this
will be stated explicitly. Similarly, when exploring a free triangulation
and a revealed triangle splits it into two components, we will specify at
each time which of the components is being further explored.

Clearly there are many degrees of freedom in carrying out the exploration
process. In an infinite map, the exploration might leave some parts
unexplored even after infinitely many steps, if the chosen boundary edges
are farther  and farther from the root. However, such pathologies are easy
to avoid, e.g.\ by guaranteeing that any boundary edge will be eventually
chosen as a base for exploration.

\medskip

We now give the probability distributions governing the exploration process
in each case. In a free triangulation of an $m$-gon, the probability that
the revealed triangle connects the edge to an internal vertex is
\[
  \frac{Z_{m+1}}{\alpha Z_m} = \frac{2m-3}{3m+3}.
\]
The probability of connecting to a vertex at distance $k$ and thus
splitting the triangulation to two free triangulations with boundary sizes
$k+1$ and $m-k$ (since there is a single common vertex) is
\[
  \frac{Z_{k+1} Z_{m-k}}{Z_m} \approx 9^{-k} Z_{k+1},
\]
where the approximation holds for fixed $k$ and $m\to\infty$ (in fact it
even holds for $k=o(\sqrt{m})$).

In the UIPT of an $m$-gon, there are similar formulas. The probability that
the revealed triangle connects the edge to an internal vertex is 
\[
  \frac{C_{m+1}}{\alpha C_m} = \frac{2m-1}{3m-3}.
\]
Finally, the probability of connecting to a vertex at distance $k$ so that
there is a free triangulation of a $k+1$-gon and a UIPT in the $m-k$-gon is
\[
  \frac{Z_{k+1} C_{m-k}}{C_m} \approx 9^{-k} Z_{k+1}.
\]
Note that the asymptotics are identical for the free triangulation and the
disc UIPT.

For the full plane UIPT, note that the whole map with a single triangle
removed is just a disc UIPT with boundary size 3. Thus we start with the
root triangle and proceed to explore the complement which is a disc UIPT
(this was done in \cite{UIPT2} in order to analyze the growth rate of the
UIPT).

\begin{proof}[Proof of \thmref{2segs}]
  To see that the weak limits in the theorem exist (and are equal) it is
  suffices to prove that for any possible finite neighborhood $B$ of the
  root edge, $\mu_m(B)$ and $\nu_m(B)$ have the same limit (abusing
  notation, $\mu_m(B)$ is the measure of all triangulations containing the 
  neighborhood $B$). That the limit is a probability measure follows using
  dominated convergence in exactly the same manner as was used for the full
  plane UIPT in \cite{UIPT1}.

  Any finite neighborhood $B$, can naturally be revealed by a finite number
  of exploration steps. The probability under $\mu_m$ or ($\nu_m$) of
  finding $B$ is just the probability that each of the exploration steps
  reveals a triangle connected to the same vertex as in $B$. Thus it
  suffices to show that for each step, the probability of finding a
  particular triangle at that step (conditioned on having succeeded so far)
  converges to some limit.
    
  This fact is easily verified from the distributions of steps of the
  exploration process, which are explicitly known. After a some number of
  steps have been made, with the revealed triangles matching those required
  for $B$, the boundary length is $m+d$ for some $d$ depending only on $B$
  and on the order in which its presence is verified. The probability for a
  disc UIPT of finding a triangle with a new vertex at this stage is
  \[
    \frac{C_{m+d+1}}{\alpha C_{m+d}} \to \frac{2}{3}.
  \]
  The probability of finding a triangle connecting to a boundary vertex at
  distance $k$ is
  \[
    \frac{Z_{k+1} C_{m+d-k}}{C_{m+d}} \to 9^{-k} Z_{k+1}.
  \]
  
  Thus we see that the limit of $\nu_m$ exists. For $\lim \mu_m$ the
  formulas are the same with $Z$'s replacing all $C$'s, and the same
  limits hold.
\end{proof}

\subsection{Markovian nature of the half plane UIPT}

A consequence of the exploration process for the half plane UIPT is the
following property:

\begin{thm}\label{T:markovian_uipt} 
  Let $T$ be a sample of the $\nu_\infty$, and let $e$ some boundary edge.
  Remove from $T$ the triangle containing $e$, and denote the infinite
  2-connected component of the remaining map by $T'$. Then $T'$ also has law
  $\nu_\infty$, and is independent of the configuration of the removed
  triangle in $T$.
\end{thm}

This fact, that the distribution of the remaining map is identical for all
steps of the exploration, makes the half plane UIPT the easiest to
analysis. Markov chains in other settings are replaced by a random walk
with fixed step distributions. The step distribution is asymmetric and has
infinite variation, but the cumulative effect of a number of steps is still
just the sum of independent random variables (in the half plane case), and
the vast body of knowledge concerning such sums can be tapped.

This property, even with the added restriction of shift invariance, does
not uniquely define the half plane UIPT. There is a one parameter family of
measures with the above property, with the half plane UIPT being an extreme
case. These measures may be defined as weak limits of finite measures,
weighted so as to have fewer internal vertices.

\begin{proof}
  Consider a step of the exploration process for a sample of $\nu_m$, where
  the triangle containing a boundary edge $e$ is revealed. As noted in the
  previous proof, the probability of finding a triangle with an internal
  vertex, or a triangle with boundary vertex at distance $k$ from $e$
  converge to some limit probabilities. The remaining map either has law
  $\nu_{m+1}$ or contains a part with law $\nu_{m-k}$. In either case, as
  $m\to\infty$, the law of the remaining component converges to
  $\nu_\infty$.
\end{proof}

\begin{figure}
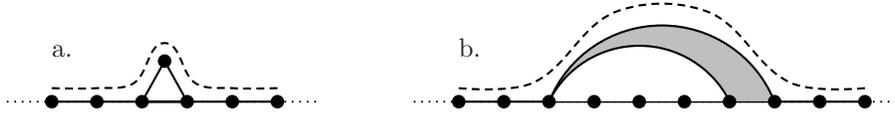

\begin{center}
  \psset{unit=.6mm, dotsep=2pt, dash=3pt 2pt}
\subfigure{
  \pspicture(0,0)(70,20)
    \put(10,10){a.}
    \psline[linestyle=dotted](0,0)(70,0)
    \psline(60,0)(30,0)(35,9)(40,0)(10,0)
    \multips{0}(10,0)(10,0){6}{\pscircle*(0,0){1.5}}
    \pscircle*(35,9){1.5}
    \pscurve[linestyle=dashed](10,3)(20,3)(28,4)(35,13)(42,4)(50,3)(60,3)
  \endpspicture
}
\hfill
\subfigure{
  \pspicture(0,0)(110,20)
    \put(10,10){b.}
    \psline[linestyle=dotted](0,0)(110,0)
    \psline(10,0)(100,0)
    \psarc[fillstyle=solid,fillcolor=lightgray](55,-10){26.9}{22}{158}
    \psarc[fillstyle=solid,fillcolor=white](50,-10){22.4}{27}{153}
    \multips{0}(10,0)(10,0){10}{\pscircle*(0,0){1.5}}
    \pscurve[linestyle=dashed](10,3)(26,4)(45,20)(65,20)(84,4)(100,3)
  \endpspicture
}
\caption{\label{fig:half_plane_exp} Exploration of the half plane UIPT: In
(b), the lower part contains an independent Boltzmann triangulation of the
polygon. In both cases, the part above the new boundary (the dotted line)
contains a half plane-UIPT, independent of the first triangle.}
\end{center}
\end{figure}

Thus the exploration process of the half plane UIPT is as follows: At each
step a boundary edge may be chosen in an arbitrary manner (which may be
either random or deterministic in terms of previous choices). The triangle
containing that edge has the following distribution: With probability 2/3
the third vertex is an internal vertex. Otherwise, with probability
\[
  p_k = 9^{-k} Z_{k+1} = \frac{(2k-2)!}{4^k (k-1)!(k+1)!}
\]
the third vertex is a boundary vertex at distance $k$ to the right (or
left) of the base edge. It is not hard to check that $2/3 + \sum 2p_k = 1$.
In the latter case, the $k+1$-gon that appears contains a free
triangulation. In either case, the remaining region with the infinite
boundary contains a half-plane UIPT that is independent of both the choice
of third vertex and of the free triangulation (see
\figref{half_plane_exp}).

\subsection{The Airy stable process}

The Airy stable process (ASP for short) is a particular instance of a
stable process, and is a L\'evy process where the independent increments
have the Airy distribution, i.e.\ are totally asymmetric stable random
variables with index 3/2. An ASP $X_t$ is scale invariant in that for
$\lambda>0$ the process $\lambda^{-1} X_{\lambda^{3/2} t}$ has the same law
as $X_t$. Sample paths are a.s.\ right continuous and have left limits
satisfying $X_{t^-} \ge X_t$ (i.e.\ the process has negative jumps but no
positive jumps). Thus ASP's have two properties that simplify their
analysis, namely being scale invariant, and spectrally negative (having
only negative jumps). For further details see \cite{zolo, bert}.

The following proposition relates the ASP to the scaling limit of random
walks:
\begin{prop}
  Let $S_n = \sum^n X_i$ be the sum of i.i.d.\ real random variables with 0
  expectation, satisfying
  \begin{align*}
    \P(X>t)  &= o(t^{-3/2}),  \\
    \P(X<-t) &= ct^{-3/2}(1+o(1))
  \end{align*}
  as $t\to\infty$ for some $c>0$, then there exists the scaling limit
  \[
    \lambda^{-1} S_{\lfloor \lambda^{3/2}t \rfloor} \to X_{kt}
  \]
  for some constant speed $k$.
\end{prop}

\subsection{Exploration and percolation interfaces}

Suppose that each vertex of a triangulation is colored black or white. This
partitions the vertices into connected monochromatic clusters ---
percolation clusters. Note that each cluster has a color, and every vertex
is part of some cluster (as opposed to the common terminology where
clusters are of only a single color.) Of course, a percolation cluster can
not be adjacent to a cluster of the same color.

Because of the local combinatorics of the maps, specifically that all faces
are triangles, whenever two percolation clusters of opposite colors are
adjacent, there is a well defined interface between them. The interface is
a path in the dual graph, that uses exactly the dual edges of edges with
differently colored endpoints. Since faces are triangles, dual vertices
have degree 3, and either two or none of the three dual edges are part of
an interface. Thus the interfaces form cycles and infinite lines, and can
not intersect each other. 

An exception is the case of triangulations of a disc. Such maps have an
external face that has degree $m$ (possibly infinite), and consequently
there may be more then two dual edges in interfaces connected to it. In the
case of a finite boundary, the colors along the boundary must flip an even
number of times. In the infinite case, the number may be any finite number,
or $\infty$. When an interface reaches the external face, we will say that
it terminates there, and the duals of bi-colored boundary edges will be
called interface end-points. Thus in the presence of a boundary, there are
two other types of interfaces: having two end-points, and having a single
end-point and being infinite (in Bernoulli percolation on a UIPT, the
latter have probability 0).

\medskip

Since a percolation interface passes from a triangle to an adjacent
triangle of the map, it is natural to explore the triangulation along a
percolation interface. Pick a base edge for an exploration step to be (the
dual of) a bi-colored boundary edge. The interface using that edge must
leave the triangle through (the dual of) one of the other two edges. If the
exploration process is used to reveal that triangle, the other two edges
become boundary edges of the remaining unknown part of the map, and so it
is always possible to continue to explore along the same interface. The
exception is of course that the interface may reach another endpoint along
the boundary, in which case the interface has been revealed completely.

\begin{figure}
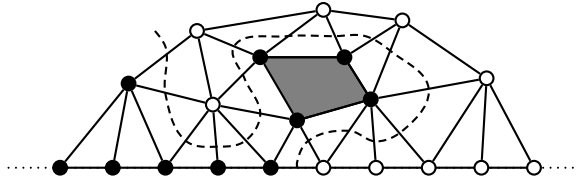
\begin{center}
  \newcommand{\verW}{\Cnode[fillstyle=solid,fillcolor=white]}
  \newcommand{\verG}{\Cnode[fillstyle=solid,fillcolor=gray]}
  \newcommand{\verB}{\Cnode[fillstyle=solid,fillcolor=black]}
  \psset{unit=.7mm, dash=3pt 2pt, dotsep=2pt, radius=1.5}
  \pspicture(-10,0)(100,30)
    \psline[linestyle=dotted](-10,0)(100,0)
    \psline(0,0)(90,0)

    \psline[fillstyle=solid,fillcolor=gray](45,9)(59,13)(54,21)(38,21)

    \verB( 0,0){b0} \verB(10,0){b1} \verB(20,0){b2} \verB(30,0){b3}
    \verB(40,0){b4} \verW(50,0){b5} \verW(60,0){b6} \verW(70,0){b7}
    \verW(80,0){b8} \verW(90,0){b9}

    \verB(13,16){v0} \verW(29,12){v1} \verB(45,9){v2} \verB(59,13){v3}
    \verW(81,17){v4} 

    \verW(26,26){v5} \verB(38,21){v6} \verB(54,21){v7} \verW(65,28){v8}
    \verW(50,30){v9}

    \ncline{b0}{v0} \ncline{b1}{v0} \ncline{b2}{v0} \ncline{b2}{v1}
    \ncline{b3}{v1} \ncline{b4}{v1} \ncline{b4}{v2} \ncline{b5}{v2}
    \ncline{b5}{v3} \ncline{b6}{v3} \ncline{b7}{v3} \ncline{b7}{v4}
    \ncline{b8}{v4} \ncline{b9}{v4}

    \ncline{v0}{v1} \ncline{v1}{v2} \ncline{v2}{v3} \ncline{v3}{v4}
    \ncline{v0}{v5} \ncline{v1}{v5} \ncline{v1}{v6} \ncline{v2}{v6}
    \ncline{v3}{v7} \ncline{v3}{v8} \ncline{v4}{v8} \ncline{v5}{v9}
    \ncline{v6}{v9} \ncline{v7}{v9} \ncline{v8}{v9} \ncline{v5}{v6}
    \ncline{v6}{v7} \ncline{v7}{v8}

    \ncline{vg}{v2} \ncline{vg}{v7} \ncline{vg}{v6}
    
    \pscurve[linestyle=dashed](45,0)(45,1)(47,5)(55,7)(60,5)(65,7)
            (70,15)(63,22)(58,25)(52,25)(43,25)(33,23)(35,17)(38,10)
            (35,5)(30,4)(24,5)(20,14)(20,23)(18,26) 
  \endpspicture

  \caption{\label{fig:exploration} The exploration process can proceed
    along the percolation interface (dotted) revealing the triangles in it.
    The gray area can no longer be visited by the interface.}
\end{center}\end{figure}

An important fact to note is that when a map is explored along percolation
interfaces, the number of interface end-points along the boundary does not
increase. When an internal vertex is encountered, the number remains the
same as the interface continues to one side or the other of the vertex.
When the map is split into two independent components, the interfaces
end-points are split up between the two boundaries. Typically we will have
a boundary that is composed of some number of segments of alternating
colors, and we will track the evolution of the segment lengths (in the
boundary of the unknown component) as the map is explored. The above means
that the number of segments does not increase.

\section{Half plane crossing probabilities} \label{sec:half_plane}
%%%%%%%%%%%%%%%%%%%%%%%%%%%%%%%%%%%%%%%%%%%

We consider probability spaces of the following type. In a half plane UIPT,
internal vertices are colored black or white with probability 1/2 each
independently of all other vertices (i.e.\ critical Bernoulli percolation).
Vertices of the boundary, consisting of an infinite path, are given
predetermined colors. 

Once the boundary colors are fixed, so is the number of mono-chromatic
segments along the boundary, and hence the number of end-points of
percolation interfaces. An interface that terminate at a given end-points
can be either infinite or end up at a second end-point. From the analysis
below it will be clear that the number of infinite percolation interfaces
is a.s.\ 0 or 1, depending on the parity of the number of end-points along
the boundary.

The original motivation for considering the percolation interfaces is the
study of crossing probabilities. The two are simply connected, since if one
knows how the interfaces connect the end-points in pairs, it is immediate
to check weather two boundary segments of the same color are in the same
percolation cluster or not. If there is an interface that separates them
them they are obviously not in the same cluster since any path from one to
the other must cross the interface, and so cannot be of a single color.
Conversely, if no interface separates the two segments, then it is easy to
find a path connecting the two to a single percolation cluster, essentially
by following parallel to interfaces, on the side having the appropriate
color.

\begin{figure}
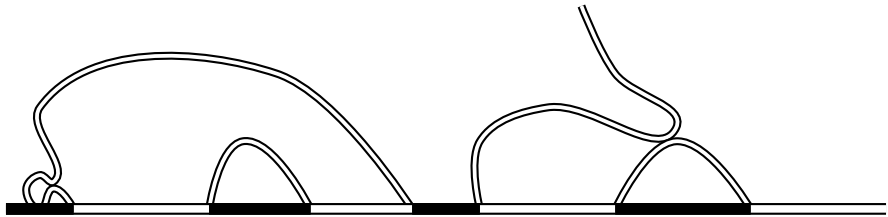

\begin{center}
  \psset{unit=.9mm, dotsep=2pt, dash=3pt 2pt}
  \pspicture(0,0)(130,30)
   {\psset{doubleline=true,doublesep=.5mm}
     \pscurve(10,0)(7,3)(5,0)(3,3)(5,5)(7,4)(5,15)(40,20)(60,0)
     \pscurve(30,0)(35,10)(45,0)
     \pscurve(70,0)(70,10)(80,15)(99,12)(90,20)(85,30)
     \pscurve(90,0)(99,10)(110,0)
   }
   {\psset{doubleline=true,doublesep=1mm}
     \psline(0,0)(130,0)
     {\psset{doublecolor=black} 
       \psline(0,0)(10,0)
       \psline(30,0)(45,0)
       \psline(60,0)(70,0)
       \psline(90,0)(110,0)
     }
   }
  \endpspicture
  \caption{\label{fig:interfaces} Each interface has a black path on one
  side and a parallel white path on the other, with one infinite interface
  (if the infinite segments are of the same color, there is a.s.\
  no infinite interface).}
\end{center}
\end{figure}

Thus in \figref{interfaces}, the two leftmost white segments are part of
the same percolation cluster, and the two rightmost white segments are part
of a second white percolation cluster (with a path passing between the two
touching interfaces). Two interfaces can not pass through the same point,
not can they pass through boundary points (since they are paths in the dual
map). When interfaces bounce off of the boundary, themselves or other
interfaces (as in the figure), they actually visit a neighboring dual
vertex and not the same vertex.

The process is a random map analogue of the percolation exploration process
on the triangular lattice in the plane \cite{Smir,Wern}. At any step, the
explored region has a boundary that is colored white on one side and black
on the other, and the exploration continues at the interface between the
two. The key difference is that here the exploration reveals not just the
colors of vertices but also the underlying random graph itself.

A second important difference is in the way vertices are removed from the
boundary. In the Euclidean lattice process there are steps at which a
segment is removed from the boundary. This happens whenever the process
hits its past or another part of the boundary. However, in the Euclidean
setting before touching itself the process is necessarily very close to its
past, so that such large changes are predictable. Alternatively, viewed
from infinity these moves are not discontinuities at all. In the random
lattice, the past gives no indication of these jumps. In the scaling limit,
it indeed turns out that the Euclidean case has predictable jumps w.r.t.\
the natural filtration, whereas the random lattice process gives rise to
ASP's which have (unpredictable) jumps.

\subsection{Basic properties}

We begin with the simplest crossing probability problem in the half planar
scenario. In a half plane UIPT, interior vertices are colored randomly
independently black or white with equal probabilities. The boundary is
colored as follows: All vertices to the left of some point are white. This
is followed by $a$ black vertices, $b$ white vertices, and the rest are
again black (see \figref{2segs}). We wish to evaluate the crossing
probability
\[
  Q_{a,b} = \P(\textrm{Both black segments are part of the same percolation
   cluster.})
\]
Here and in the rest of this section, $\P$ is the measure on colored half
planar UIPT's as above, with boundary colors that will be specified as
needed, and clear from the context.

Since the event of having a black crossing is increasing with respect to
making additional vertices black, it is clear that $Q_{a,b}$ is monotone
increasing in $a$ and decreasing in $b$. Also, for any given map, the event
that the black segments are connected and the event that the white segments
are connected (by a white path) are disjoint events. Additionally, as will
be shown below and similarly to the case of the Euclidean lattice, the
probability that neither color has crossing (related to the four arm
exponent) is 0. Thus up to measure 0, the two events are complements. Since
the UIPT is symmetric with respect to reflection, this implies the identity
$Q_{a,b}+Q_{b,a}=1$. A natural scaling problem is to evaluate the limit of
$Q_{\lambda a,\lambda b}$ as $\lambda\to\infty$, or at least to show the
limit exists.

When the boundary vertices are colored in this way, there are three
end-points of percolation interfaces at the boundary. Since there is
evidently an infinite percolation cluster of each color, with an infinite
interface between them, one of the end-points must belong to an infinite
interface. The other two are a.s.\ the two ends of a finite interface. This
leaves two possibilities: Either the infinite interface ends at the
rightmost end-point, or at the left-most (it a.s.\ does not end at the
center end-point, since percolation interfaces cannot intersect). By
following the finite interface on the appropriate side, one can find either
a white or a black crossing, depending on which two of the three end-points
are connected by it.

\medskip

As in the Euclidean case, the problem of determining crossing probabilities
has several forms. One important form of the problem is as follows:
Consider a (topological) triangle with vertices colored black along one
edge, white along another, and randomly along the third. The percolation
cluster of the black edge contains some of the vertices of the remaining
edge of the triangle, and the problem is to find the distribution of the
farthest of these along the third edge. A famous form of Cardy's formula
for critical percolation in the Euclidean plane is that the scaling limit
of this distribution is conformally invariant, and that if the triangle is
equilateral it is the uniform distribution on the third edge
\cite{Smir,Wern}. 

There is a similar version of finding $Q_{a,b}$ in the half plane UIPT.
Assume the boundary contain $a$ black vertices in a single segment, with
all other vertices being white. Consider the black percolation cluster
containing the black segment. It is adjacent to some boundary vertices to
its right (a.s.\ finitely many). Let $X$ be the distance along the boundary
from the black segment to the farthest of those. It is clear that
$Q_{a,b}=\P(X>b)$ (the two probability spaces can be coupled so that the
map is identical and only the boundary colors differ).

Yet another form of the same problem is to color $b$ vertices white, and
all vertices to their right black, and consider the cluster $B$ containing
the infinite black segment. This cluster is a.s.\ adjacent to infinitely
many vertices to the left of the white segment. $Q_{a,b}$ is also the
probability that the nearest of those is at distance at most $a$ from the
leftmost white boundary vertex. Note that in this form some part of the
boundary was not assigned a color. If those vertices are colored randomly,
there would be a distinction here between the cluster being adjacent to
some vertex and containing the vertex. As an aside, it is not difficult to
prove that asymptotically the two options give the same probabilities.

\subsection{A Markov chain}

\begin{figure}
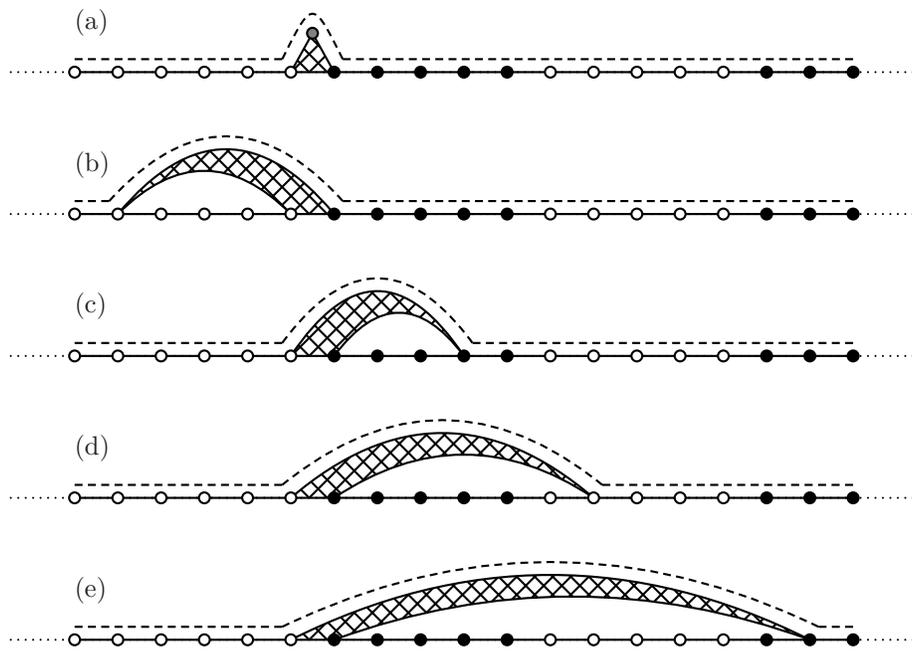

\begin{center}
  \psset{unit=.575mm, radius=1.5, dash=3pt 2pt, dotsep=2pt, fillstyle=solid}

  \newcommand{\basesegs}{ 
    \psline[linestyle=dotted](-15,0)(195,0)
    \psline(0,0)(180,0) 
    \multips{0}(  0,0)(10,0){6}{\pscircle[fillcolor=white](0,0){1.5}}
    \multips{0}( 60,0)(10,0){5}{\pscircle[fillcolor=black](0,0){1.5}}
    \multips{0}(110,0)(10,0){5}{\pscircle[fillcolor=white](0,0){1.5}}
    \multips{0}(160,0)(10,0){3}{\pscircle[fillcolor=black](0,0){1.5}} 
  }

  \subfigure{ 
    \pspicture(-15,0)(195,20)
      \psline[linestyle=dashed](0,3)(48,3)
      \psbezier[linestyle=dashed](48,3)(55,17)(55,17)(62,3)
      \psline[linestyle=dashed](62,3)(180,3)
      \psline[fillstyle=crosshatch](50,0)(55,9)(60,0)
      \pscircle[fillcolor=gray](55,9){1.5} \basesegs \put(0,10){(a)}
    \endpspicture
  }

  \subfigure{
    \pspicture(-15,0)(195,20)
      \psline[linestyle=dashed](0,3)(8,3)
      \parabola[linestyle=dashed](8,3)(35,18)
      \psline[linestyle=dashed](62,3)(180,3)
      \parabola[fillstyle=crosshatch](10,0)(35,15)
      \parabola(10,0)(30,10)
      \basesegs
      \put(0,10){(b)}
    \endpspicture
  }

  \subfigure{
    \pspicture(-15,0)(195,20)
      \psline[linestyle=dashed](0,3)(48,3)
      \parabola[linestyle=dashed](48,3)(70,18)
      \psline[linestyle=dashed](92,3)(180,3)
      \parabola[fillstyle=crosshatch](90,0)(70,15)
      \parabola(90,0)(75,10)
      \basesegs
      \put(0,10){(c)}
    \endpspicture
  }

  \subfigure{
    \pspicture(-15,0)(195,20)
      \psline[linestyle=dashed](0,3)(48,3)
      \parabola[linestyle=dashed](48,3)(85,18)
      \psline[linestyle=dashed](122,3)(180,3)
      \parabola[fillstyle=crosshatch](120,0)(85,15)
      \parabola(120,0)(90,10)
      \basesegs
      \put(0,10){(d)}
    \endpspicture
  }

  \subfigure{
    \pspicture(-15,0)(195,20)
      \psline[linestyle=dashed](0,3)(48,3)
      \parabola[linestyle=dashed](48,3)(110,18)
      \psline[linestyle=dashed](172,3)(180,3)
      \parabola[fillstyle=crosshatch](170,0)(110,15)
      \parabola(170,0)(115,10)
      \basesegs
      \put(0,10){(e)}
    \endpspicture
  }

  \caption{\label{fig:2segs} The various possibilities when a triangle is
  revealed. In (a),(b),(c), the process continues with the new boundary
  along the dashed line). In (d),(e) a crossing of one of the colors has
  been revealed.}
\end{center}
\end{figure}

We return to the first variation above, where the boundary has two finite
segments of lengths $a,b$. Denote by $A$ the event that the two black
segments (of lengths $a,\infty$) are part of a single percolation cluster.
The exploration process is now used to explore the map as well as the color
of each encountered vertex, and thereby see whether $A$ occurs or not. We
may explore the map at any edge of the boundary, but our first choice is to
reveal at each step the triangle containing the edge that connects the
finite black segment to the infinite white segment to its left (as in
\figref{2segs}), and thus we explore the map along the percolation
interface ending there.

With probability 2/3 The explored triangle reveals a new vertex, which is
randomly colored black or white and inserted into the boundary where the
edge was (\figref{2segs} (a)). Looking at the new boundary's colors,
finding a white vertex leaves the boundary coloring exactly as it was.
Finding a black vertex increases by 1 the length of the black segment. Note
this introduces black vertices in the new boundary that are not part of the
original segment, but are in the same connected black cluster. Thus a
connection between these vertices and the infinite black segment
establishes the occurrence of $A$.

Another possibility is that for some $k>0$, with probability $p_k$
(defined above), the third vertex of the triangle is the white vertex at
distance $k$ to the left of the base edge. In this case, the resulting
$(k+1)$-gon, contains some finite unknown triangulation. However, since the
boundary of this $k+1$-gon contains only white vertices, this triangulation
and its coloring have no influence on $A$ (\figref{2segs} (b)). Changing
any part of the map inside the $k+1$-gon does not alter the connectivity of
percolation clusters outside the polygon. Looking at the interfaces, the
followed interface has been followed to a point where it is adjacent to the
infinite white cluster. Since none of the interfaces can enter the enclosed
region, it is inconsequential. The remaining boundary has $k$ fewer white
vertices but since these are removed from an infinite segment, the color
pattern is unchanged. Together with the 1/3 probability of adding a white
vertex this gives to a total probability of 1/2 that the boundary color
pattern remains the same.

Finally, also with probability $p_k$ the third vertex is at distance $k$ to
the right of the edge. Let $a'$ be the length of the black segment just
before that step. If $k<a'$ then just as in the previous case, the map in
the $k+1$-gon is independent of $A$. In this case the process can continue,
with a the new black segment having length $a'-k$. However, if $a'\le k <
a'+b$ then the third vertex is a white, and so the revealed triangle
contains an edge between the two white segments, hence the two white
segments are connected and $A$ does not occur. Finally, if $k\ge a'+b$ then
the exploration has found an edge between the black segments and so $A$
does occurs (see \figref{2segs} (c,d,e)). Thus $A$ is determined at the
first step of this type at which $k\ge a'$.

In terms of the percolation interfaces, in the case $k\ge a'+b$ the
interface is adjacent to the infinite black segment, and continues in the
infinite unknown component, whereas the two other endpoints are in the
boundary of the $k+1$-gon, and so must belong to the same interface. In the
case $a'\le k < a'+b$, the followed interface continues into the free
triangulation in the $k+1$-gon, and so it must connect eventually to the
central end-point that is also in the boundary of the $k+1$-gon. The
rightmost interface end-point is the only one in the infinite component,
and so must belong to the infinite interface.

\medskip

In summary, apart from the length of the black segment, the boundary color
pattern remains fixed throughout the process, until some step at which the
occurrence of $A$ is determined. The next step is to study the evolution of
the segment's length. Let $S_n$ be the length after $n$ steps of the
process have been made. The above discussion amounts to the statement that
the law of $S_n$ is that of a random walk on $\Z$ with $S_0=a$ and i.i.d.\
steps $X_i$. The distribution of the steps $X_i$ is given by
\begin{equation} \label{eq:step_dist}
  X_i = \begin{cases}
     0  & \text{with pr. $p_0 = 1/2$,} \\
     1  & \text{with pr. $p_{-1} = 1/3$,} \\
    -k  & \text{with pr. $p_k = Z_k 9^{-k} =
                              \frac{(2k-2)!}{4^k(k-1)!(k+1)!}$.} \\
  \end{cases}
\end{equation}

Occurrence of the event $A$ is determined at the hitting time of the
non-positive integers ($\Z^-$) by the random walk, denoted as $T_- = \inf
\{n\vert S_n\le 0\}$. Note that the steps have $\E X_i=0$, and therefore
a.s.\ $T_- < \infty$, and consequently there is a.s.\ a crossing in one of
the two colors (or equivalently, critical percolation dies). Thus with the
above notation, we can write $A = \{|S_{T_-}| \ge b\}$. Letting $\P_x$
denotes the probability measure for a process started at $x$, we have
proved

\begin{thm} \label{T:2segs}
  Let $\P_a$ be the probability measure of a random walk started at $a$
  with step distribution \eqref{eq:step_dist}. The crossing probability is
  given in terms of the hitting measure on $\Z^-$ by $Q_{a,b} =
  \P_a\left(|S_{T_-}| \ge b\right)$.
\end{thm}

Translating to the hitting probability variation of the problem, we find
that the probability that the percolation cluster of a black boundary
segment of length $a$ is adjacent to a boundary vertex at distance greater
than $b$ to the segment's right is also given by $\P_a\left(|S_{T_-}| \ge
b\right)$. 

As an aside, take the case where $a$ boundary vertices are colored black
and the rest are colored randomly. The color of such boundary vertices is
treated as unknown until the exploration process encounters them. As the
map is explored, the length of the black segment is described by a random
walk as before. The difference is that once the random walk becomes
non-positive it is not always killed. With probability half the hit vertex
is found to be black as well, and the process continues with a segment of
length 1 and a new and closer goal of reaching distance $b-|S_n|$. However,
starting at 1, the random walk is likely to become non-positive again
quickly, so asymptotically (as $a,b$ are scaled to infinity) this makes no
difference. The details are left to the reader.

\medskip

The UIPT may also be explored from other edges along the boundary. The
previous choice of the edge connecting the finite black segment to the
infinite white segment was convenient but arbitrary. One other choice is
equally convenient, namely, to add the triangle containing the edge that
connects the finite white segment to the infinite black segment, and
explore the map along the percolation interface terminating there. The
resulting analysis is exactly as above, except that the colors are
exchanged, as are the roles of left and right. We get the additional
identity
\[
  Q_{a,b} = \P_b\left(|S_{T_-}| < a\right),
\]
and since $Q_{a,b}+Q_{b,a}=1$, we can combine this with the previous
formula for $Q_{a,b}$ to deduce an identity concerning the hitting
probabilities of the random walk:

\begin{coro}
  \[
    \P_a\left(|S_{T_-}| < b\right) + \P_b\left(|S_{T_-}| < a\right) = 1.
  \]
\end{coro}

In this way the symmetry of the UIPT gave rise to information on the random
walk. It is interesting to note, that this identity is not a consequence of
some path transformation, neither at the discrete version, nor at the
scaled version for ASP's that will follow. This is not surprising: the
random walk describes the interface between a white percolation cluster and
a black one along which the map is explored. The two exploration processes
that give rise to the two formulas explore two disjoint interfaces. On the
event $A$, there is a black cluster spanning the finite and the infinite
black segments. This cluster has an interface with the finite white cluster
is surrounds, and with the infinite white cluster. The two exploration
processes explore these two boundaries (see \figref{interfaces}). A
similarly picture holds when $A$ does not occur.

\subsection{A Scaling limit}

Consider the crossing probability $Q_{a,b}$ when $a,b$ are large. Since the
random walk $S_n$ has a scaling limit, and the hitting measure of the limit
on $\R^-$ is non-trivial, the only limit that makes sense is to scale $a,b$
by an equal factor $\lambda$. Letting $\lambda\to\infty$, the normalized
process $\lambda^{-1}S_{\lambda^{3/2}t}$ scales to an ASP, and therefore we
get: 

\begin{thm} \label{T:sclim2}
  The scaling limit of the crossing probabilities $Q_{\lambda a,\lambda b}$
  exists and is expressed in terms of the hitting measure on $\R^-$ of an
  ASP started at $a$ by the relation
  \[
    \lim_{\lambda\to\infty} Q_{\lambda a,\lambda b}
    = \P_a\left(|Y_{T_-}| > b\right).
  \]
\end{thm}

Since ASP's are scale invariant, the limit is a function of the ratio $a/b$
(as a scaling limit must be). Indeed (see \cite{bert} for definitions and
details), since the ladder height process of $-Y_t$ is a stable
subordinator of index 1/2, it is possible to derive an explicit formula for
the hitting distribution:
\begin{equation}\label{eq:overshoot}
  \P_a\left(|Y_{T_-}| > b\right) = \pi^{-1} \arccos \frac{b-a}{a+b} .
\end{equation}

The symmetry in the roles of $a,b$ noted above in the discrete case,
translates to the identity

\begin{coro} \label{C:symmetry}
  \begin{equation}\label{eq:symmetry}
    \P_a\left(|Y_{T_-}| > b\right) = \P_b\left(|Y_{T_-}| < a\right)
  \end{equation}
\end{coro}

While this identity is easily verified from \eqref{eq:overshoot}, it would
be interesting to find a more direct proof. As noted, a proof using simple
path transforms seems unlikely: \eqref{eq:symmetry} does not hold for
arbitrary L\'evy processes. Moreover, the LHS event is correlated to having
large negative jumps while the RHS even is negatively correlated to such
jumps. Among stable processes, \eqref{eq:symmetry} holds for (and only for)
those processes $X$ s.t.\ the ladder height process of $-X$ is a stable
subordinator of index 1/2. These processes are stable processes with index
$\alpha \in [1/2,3/2]$ and suitable asymmetry parameter $\beta =
\frac{\tan(\alpha-1)\pi/2}{\tan\alpha\pi/2}$.

\subsection{A third approach}

A third interesting way to assess the probability $Q_{a,b}$ is to explore
the UIPT along the percolation interface that starts at the edge between
the two finite boundary segments. We then need to determine which of the
two other end-points the interface reaches. As long as it encounters only
the finite segments, the interface will bounce outside towards the infinite
map and there is no topological change. The moment a triangle is revealed
that forces the interface to pass by one of the infinite segments, the
interface will continue in the resulting free triangulation, and terminate
at the second end-point in the boundary of the finite map. Thus determining
occurrence of $A$ is equivalent to finding out which of the two infinite
segments is hit first by the interface.

Considering the lengths of boundary segments, with this approach there are
no steps that effect only an infinite segment and make no difference. Each
step of the exploration process along this interface effects the length of
exactly one of the two finite segments, with probability 1/2 each. Each
length increases by 1 with probability 1/3 and decreases by $k$ with
probability $p_k$. If in this way one of the lengths becomes negative ---
i.e.\ a triangle is revealed that connects the edge to some boundary vertex
far enough to be in either of the infinite segments --- then a crossing has
been established in one of the two colors, and occurrence of $A$ has been
determined.

Let $S_n,S'_n$ denote the segment lengths after $n$ steps, so that $S_0=a$
and $S'_0=b$. As before, $S_n$ and $S'_n$ are random walks on $\Z$ with
step distribution given by \eqref{eq:step_dist}. The two walks are not
independent since at each step of the exploration process, exactly one of
them makes a non-zero step. In terms of the percolation interface, finding
which color has a crossing is equivalent to finding which infinite segment
is touched first by the process. If $T_-$ and $T'_-$ are the corresponding
hitting times of $\Z^-$ by the two random walks, then we have $A = \{T_- >
T'_-\}$ (noting that equality is impossible).

The dependence between the two random walks is only a minor inconvenience.
Instead of discrete time let the exploration proceed with a Poisson clock,
so that triangles are revealed at rate 1. Since the two random walks do not
change their value at the same step, in continuous time they become
independent random walks in continuous time, progressing at rate 1 each
with the same step distribution as before. The two random walks are started
at $a,b$ and $A$ occurs if the latter hits $\Z^-$ before the former does.

\medskip

In the scaling limit, this approach simplifies further. The Poisson process
governing the progress of each random walk scales to linear time, and so
has no effect on the scaling limit. The two random walks scale to
independent ASPs, $Y,Y'$ started at $a,b$ respectively.

\begin{thm} \label{T:sclim2a}
  Let $\P_{a,b}$ be the probability measure for two independent ASP's
  started at $a$ and $b$. Let $T_-$ and $T'_-$ denote the respective
  hitting times of $\R^-$ by the ASP's, then
  \[
    \lim_{\lambda\to\infty} Q_{\lambda a,\lambda b} = \P_{a,b}(T_->T'_-).
  \]
\end{thm}

Comparing this with \thmref{sclim2} yields a second non trivial statement
concerning ASPs. Because of the scale invariance of the ASP, $T_-$ is
distributed as $a^{3/2}\tau$ where $\tau$ is the hitting time of $\R^-$ by
an ASP started at 1. Combining this with \eqref{eq:overshoot}, we find a
remarkably simple formula for the distribution of $\tau/\tau'$:

\begin{coro} \label{C:id2}
  If $\tau$ is hitting time of $\R^-$ by an ASP started at 1, and
  $\tau,\tau'$ are i.i.d., then
  \[
    \P\left(\frac{\tau}{\tau'} > (\frac{a}{b})^{3/2} \right) =
    \pi^{-1}\arccos \frac{a-b}{b+a} ,
  \]
  and hence
  \[
    \P\left(\frac{\tau}{\tau'} > t \right) =
    \pi^{-1}\arccos \frac{t^{2/3}-1}{t^{2/3}+1} .
  \]
\end{coro}

This identity is even more specific to ASP's then the previous on. Here, it
is not enough that the ladder height process be a stable subordinator of
index 1/2 to ensure this identity holds. Note though that in general the
distribution of $X/X'$ does not determine the distribution of a random
variable $X$ --- not even up to a multiplicative constant.

\subsection{Mixed growth}

When a map is being explored, the edge at which a triangle is revealed does
not have to be a deterministic function of the known part of the map. It is
easier to explore a map along percolation interfaces. We have seen above
the consequences of exploring an interface starting at each of the three
end-points along the boundary. However, it is also possible to have a mixed
exploration, along two of the interfaces or all three.

Using Poisson timing as before, assume that triangle are revealed at rate 1
at each of the two edges connecting an infinite segment to a finite one. As
before, the lengths $S_t,S'_t$ of the black and white finite segments at
time $t$ form a pair of independent random walks in continuous time, with
step distribution \eqref{eq:step_dist}, and existence of a black crossing
is determines at time $\min(T_-,T'_-)$. However, the condition for having a
crossing is different. There are two walks that can hit $\Z^-$, and two
ways a black crossing can be found. There is a black crossing either if
$T_- < T'_-$ and $|S_{T_-}| \ge S'_{T_-}$, or if $T'_- < T_-$ and
$|S'_{T'_-}| < S_{T'_-}$.

This process too has a scaling limit, giving yet another formula for the
scaling limit of $Q_{a,b}$:

\begin{coro} \label{C:id3}
  For two independent ASP's $Y,Y'$, hitting $\R^-$ at times $T_-$ and
  $T'_-$, if $\tau = \min(T_-,T'_-)$, then
  \begin{multline*}
      \P_{a,b} \big( \{\tau = T_-\}  \cap \{Y_\tau + Y'_\tau < 0\} \big)  \\
    + \P_{a,b} \big( \{\tau = T'_-\} \cap \{Y_\tau + Y'_\tau > 0\} \big)
    = \pi^{-1} \arccos \frac{b-a}{b+a} .
  \end{multline*}
\end{coro}

However, \corref{id3} is not hard to deduce directly from previous results:
Subtracting this from the identity of \thmref{sclim2a}, we get
\[
  \P \big( \{\tau = T_-\}  \cap \{Y_\tau + Y'_\tau < 0\} \big) =
  \P \big( \{\tau = T'_-\} \cap \{Y_\tau + Y'_\tau < 0\} \big).
\]
This identity can be proved directly (easy), and unlike the previous ones
holds for any two independent instances of any Levy process, since from any
state of the pair of processes, the two events occur at the same rate. Thus
\corref{id3} follows from \thmref{sclim2a}.

The choice of exploring at rate 1 at each of the two external edges can be
generalized. If triangles are revealed at different rates at the two edges,
then the two random walks will progress at corresponding rates. In the
scaling limit, the random walks converge to ASPs with distinct time rates.

Even more generally, it is possible for the exploration to proceed at a
changing rate that is dependent on the past of the processes. If the choice
of rate is such that it has a scaling limit, then an identity for the limit
results. For example one process might progress at a linearly increasing
rate, while the other might advance at a rate proportional to $|Y_t-Y'_t|$.
Thus the conclusion of \corref{id3} holds under rather weak conditions on
the processes $Y,Y'$. Eschewing formal definitions, this may be stated as:

\begin{coro}\label{C:id4}
  For a pair of ASP's $X_s,X'_s$, let $Y_t=X_{s(t)}$ (resp. $Y'_t =
  X'_{s'(t)}$) be a time changed process, with time changes $s(t),s'(t)$
  adapted to the natural common filtration of $Y,Y'$. For $T_-,T'_-,\tau$
  as before,
  \begin{multline*}
      \P_{a,b} \big( \{\tau = T_-\}  \cap \{Y_\tau + Y'_\tau < 0\} \big)  \\
    + \P_{a,b} \big( \{\tau = T'_-\} \cap \{Y_\tau + Y'_\tau > 0\} \big)
    = \pi^{-1} \arccos \frac{b-a}{b+a} .
  \end{multline*}
\end{coro}

Finally, if the edge connecting the two finite segments is used as well,
the rate of each of the random walks is increased by some equal amount (as
is the rate of the ASPs in the limit). Each random walk can become
non-positive either as a result of a contribution from the common edge or
from its unique edge. The occurrence of $A$ depends on the overshoot as well
as on which side the contribution that made the random walk non-positive
came from. In general, this leads to convex combinations of the identities
of Corollaries~\ref{C:id2} and \ref{C:id3}.

\section{Crossing probabilities in Boltzmann maps}  \label{sec:boltzmann}
%%%%%%%%%%%%%%%%%%%%%%%%%%%%%%%%%%%%%%%%%%%%%%%%%%

\subsection{An auxiliary problem} \label{ssec:3segs}

In general, planar crossing probabilities are defined for topological
rectangles in the plane, i.e.\ a simply connected planar domain, where the
(topological) boundary is separated into four segments forming the sides.
There is no real difficulty with having part of the boundary being
infinite, in which case the point at infinity may be of the corners of the
rectangle (as above), or may be in the interior of one of the sides. It
turns out that in order to analyze percolation crossing probabilities in a
Boltzmann triangulation, it is helpful first to consider crossings in the
half planar case when infinity is not one of the rectangle's corners.

To this end, we consider the case where in a half plane UIPT, there are two
finite black segments along the boundary, with all other vertices being
white (compare with the previous case where one of the segments was
infinite). There is also of course the dual case of two finite white
segments, with a finite and two infinite black segments, where crossing
probabilities are the complement of the discussed case. 

Assume then, that in the boundary there are two black segments of lengths
$a,c$ with exactly $b$ white vertices separating them, and that all other
boundary vertices are white. Internal vertices of the UIPT are colored
white or black with equal probability. A.s.\ there is a unique infinite
white cluster and no infinite black cluster. As before, let $A$ be the
event that there is a black crossing, i.e.\ the finite black segments are
part of a single percolation cluster. The complementary event (up to
measure 0) is that the finite white boundary segment is part of the
infinite white percolation cluster. Denote $Q_{a,b,c} = \P(A)$.

There are alternative definitions of $Q_{a,b,c}$. One equivalent definition
is that $Q_{a,b,c}$ is the probability that the percolation cluster of a
black boundary segment of length $a$ is adjacent to a boundary vertex to
its right at a distance in the interval $(b,b+c]$. Yet another way involves
the lack of a white crossing: Given a white boundary segment of length $b$,
that the percolation cluster containing it is not adjacent to any boundary
vertex to its left at distance greater than $c$, nor to any vertex to its
right at distance more than $a$. This two definitions are similar to the
unlimited, one sided definition for $Q_{a,b}$. Of course, there is a
symmetry between the roles of $a$ and $c$ here, so $Q_{a,b,c}=Q_{c,b,a}$
can be defined in terms of the cluster of a segment of length $c$ being
adjacent to some vertex to its left at distance in $(b,a+b]$.

The percolation interface representation is also nice. There are four
end-points which must be connected in pairs to form two finite percolation
interfaces. There are two ways to do this since the interfaces cannon
intersect, and the two correspond to $A$ and its complement.

\medskip

Let us explore the map along one of the percolation interfaces. There are
four end-points to choose from, and our first try is to explore along the
interface between the infinite white segment on the left and the black
segment of length $a$. As the map is explored along this interface, the
length of the black segment evolves as a random walk $S_n$ with step
distribution \eqref{eq:step_dist}. Suppose $S_n$ hits $\Z^-$ at time $T_-$.
If $|S_{T_-}| \in [0,b)$ then the process has revealed a white crossing,
and $A$ cannot occur. Similarly, if $|S_{T_-}| \in [b,b+c)$, then a black
crossing has been found and $A$ occurs. The limit probabilities of these
two cases when $a,b,c$ are scaled have a simple and closed form in terms of
the overshoot of an ASP (the limits are $\P_a\left(|Y_{T_-}| < b \right)$
and $\P_a\left(|Y_{T_-}| \in [b,b+c] \right)$ respectively).

However, a new difficulty arises in the event that $|S_{T_-}| \ge b+c$. If
such a triangle is encountered, then the remaining infinite unexplored
component of the UIPT has a completely white boundary, and is independent
of $A$. Existence of a crossing is now determined by the combinatorics and
colors of the finite free triangulation that is in the resulting polygon.
It is possible to continue the exploration process in a free triangulation,
and this approach is also mentioned below. The same difficulty arises when
exploring along the interface next to the second infinite white segment.

\begin{thm} \label{T:3segs}
  Let $S,S'$ be two independent ASPs hitting $\Z^-$ at times $T_-,T'_-$
  respectively, with $\tau = \min(T_-,T'_-)$, then
  \begin{align*}
    Q_{a,b,c}
      &= \P_{a,b} \big( \{\tau = T'_-\} \cap \{|S'_\tau| < c\} \big) \\
      &= \P_{c,b} \big( \{\tau = T'_-\} \cap \{|S'_\tau| < a\} \big) .
  \end{align*}
\end{thm}

\begin{proof}
  To circumvent the above difficulty, we explore along the interface
  starting at one of the other two end-points, namely at one of the edges
  that connect two finite boundary segments. Suppose the edge between the
  segments with initial lengths $a,b$ is used. As above, the lengths of
  these two segments now form two independent random walks $S_t,S'_t$ in
  continuous time, started at $a$ and $b$ respectively. If they hit $\Z^-$
  at times $T_-$ and $T'_-$ respectively, then we can express the event $A$
  as 
  \[
    A = \big\{ T'_-<T_- \big\} \cap \big\{ |S'_{T'_-}| < c \big\} .
  \]
  To get the second identity, explore an interface starting with the other
  end-point --- between the segments of length $b,c$. This simply exchanges
  the roles of $a$ and $c$, i.e.\ $S_0=c$ instead of $a$, and $a$ replaces
  $c$ in the condition for $A$ to occur.
\end{proof}

By taking a scaling limit of the above formulas, we get

\begin{coro} \label{C:3segs}
  Let $Y,Y'$ be two independent ASPs hitting $R^-$ at times $T_-,T'_-$
  respectively, with $\tau = \min(T_-,T'_-)$, then
  \begin{align*}
    \lim_{\lambda\to\infty} Q_{\lambda a,\lambda b,\lambda c}
      &= \P_{a,b} \big( \{\tau = T'_-\} \cap \{|Y'_\tau| < c\} \big) \\
      &= \P_{c,b} \big( \{\tau = T'_-\} \cap \{|Y'_\tau| < a\} \big) .
  \end{align*}
\end{coro}
  
We do not know of a closed form for this probability, nor of a direct proof
that the two formulas are equal that does not involve the natural symmetry
of the UIPT (though, conceivably, both a closed form and a proof are within
reach).

\subsection{Boltzmann triangulations}

In order to extend the analysis of crossing probabilities to percolation on
Boltzmann triangulations, we now use the fundamental relation between these
and the UIPT, specifically in the half plane. If the UIPT is conditioned on
containing some sub-triangulation with several unknown components, then the
finite components are all independent and have the Boltzmann distribution.
This allows us to express the probability of events in free triangulations
in terms of probabilities of events in the half plane UIPT.

As an example, consider the measure $\tP$ for a free triangulation in an
$m$-gon, and the event $A$ that the triangle containing a specified
boundary edge connects that edge to an internal vertex (rather than to some
other boundary vertex). The formulas governing the exploration process of a
free triangulation tell us that $\tP(A) = \frac{Z_{m+1}}{\alpha Z_m} =
\frac {2m-1}{2m-3}$. We now proceed to calculate this probability using
only exploration of the half plane UIPT.

Consider the half plane UIPT, and pick some edge $e$. Let the event $B$ be
the event that the triangle incident on $e$ connects to a vertex at
distance exactly $m-1$ to the right from $e$. From the formulas for the
half plane exploration process, $\P(B) = p_{m-1} = 9^{1-m}Z_m$. Let $e'$ be
the edge to the right of $e$, and let the event $C$ be the event that the
triangle incident on $e'$ connects it to an internal vertex of the
triangulation, so that $\P(C) = p_{-1} = 2/3$.

Now, conditioned on $B$, the finite component containing $e'$ is just a
free triangulation of a $m$-gon, and therefore
\[
  \tP(A) = \P(C|B) = \frac{\P(B|C)\P(C)}{\P(B)}.
\]
Now, $\P(B)$ and $\P(C)$ are known. To find $\P(B|C)$, note that
conditioned on $C$, removing the triangle incident on $e'$ leaves a half
plane UIPT, and the target vertex of $B$ is now at distance $m$ from $e$.
Thus $\P(B|C) = p_m = 9^{-m}Z_{m+1}$, and we conclude that
\[
  \tP(A) = \frac{9^{-m}Z_{m+1} \cdot 2/3}{9^{1-m} Z_m}
         = \frac{Z_{m+1}}{(27/2) Z_m},
\]
in agreement with the exploration in the Boltzmann map. This example is
perhaps the simplest event in a free triangulation. It involves a single
triangle, and does not involve the percolation on the map at all, but only
the map itself. The same technique extends to more complex events with
little modification.

\medskip

Turning back to crossing exponents, consider an $m$-gon whose boundary
is composed of four segments of lengths $a,b,c,d$ colored alternately black
and white, with the segments of lengths $a,c$ being black. Let $\tP$ denote
the probability measure for the Boltzmann triangulation of this disc, with
internal vertices supporting critical percolation. Let $\tA$ be the event
that the two black segments are part of the same percolation cluster. Our
goal is to find the crossing probability
\[
  \tQ = \tQ_{a,b,c,d} =  \tP(\tA).
\]

As before, let $S,S'$ be two random walks in continuous time with step
distribution \eqref{eq:step_dist}, and let $T_-,T'_-$ be their hitting
times of $\Z^-$ and $\tau = \min(T_-,T'_-)$. Denote by $\E_{a,b}$
expectation for the random walks $S,S'$ started at $a,b$ respectively.

\begin{thm}\label{T:boltzmann}
  With the above notations, the crossing probability in a Boltzmann
  triangulation is
  \[
    \tQ_{a,b,c,d} 
    = \E_{a,b} \left( \frac{p_l}{p_{m-1}}, 
                      \big\{ \tau = T'_- \big\} \cap
                      \big\{ |S'_\tau| < c \big\} \right),
  \]
  (i.e.\ expectation restricted to the event) where $l = S_\tau + S'_\tau
  +c+d-1$. 
\end{thm}

Like the crossing probabilities for the half planar case, these crossing
probabilities also have a scaling limit. Let $tQ_\lambda = tQ_{\lambda a,
\lambda b, \lambda c, \lambda d}$ be the probability of a black crossing
in the free triangulation of a disc, with boundary segments of lengths
$\lambda a,\lambda b,\lambda c,\lambda d$ (rounded to integers). Since
$p_l\approx cl^{-5/2}$, factoring out $c\lambda^{-5/2}$ from $p_l$ and
$p_{m-1}$ gives

\begin{coro} 
  The scaling limit of the crossing probabilities in Boltzmann
  triangulations is 
  \[
    \lim_{\lambda\to\infty} \tQ_\lambda = \E_{a,b} \left(
    \left( \frac{Y_\tau + Y'_\tau +c+d}{a+b+c+d} \right)^{-5/2}, 
               \big\{ \tau = T'_- \big\} \cap
               \big\{ |T'_\tau| < c \big\}  \right),
  \]
  where $\E_{a,b}$ denotes expectation for two independent ASP's started
  at $a,b$, and $T_-,T'_-,\tau$ are as before.
\end{coro}

\begin{proof}[Proof of \thmref{boltzmann}]
  We consider again the measure $\P$ for a half planar UIPT, with critical
  percolation on internal vertices, and with the boundary coloring of the
  previous sub-section, having three finite segments. Thus we have two
  black segments, the left one of length $a$ and the right one of length
  $c$, separated by $b$ white vertices and with all other boundary vertices
  being white as well.

  Denote by $e_0$ the edge connecting the segment of length $a$ to the
  infinite white segment, and by $v_0$ the boundary vertex at distance
  $m-1$ to the right of $e_0$ (where $m=a+b+c+d$). Define the auxiliary
  event $B$ as the event that the triangle supported on $e_0$ has $v_0$ as
  its third vertex. Conditioned on $B$, the interval from $e_0$ to $v_0$
  together with an edge of the triangle connecting $e_0$ to $v_0$ surrounds
  a finite component of the UIPT, which is therefore distributed as a
  Boltzmann triangulation of an $m$-gon (see \figref{boltzmann}). Since the
  boundary colors agree as well, conditioned on $B$, this component has law
  $\tP$.

  \begin{figure}
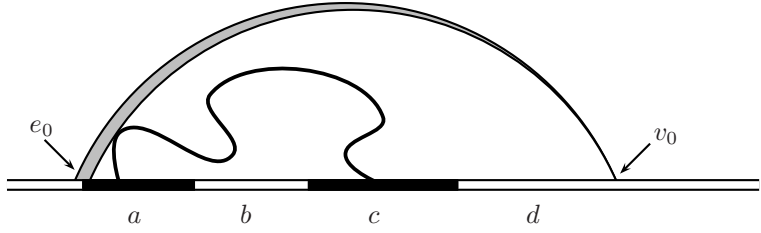

  \begin{center}
    \psset{unit=1mm, dotsep=2pt, dash=3pt 2pt, doublesep=1mm}
    \pspicture(0,-5)(100,25)
      \psarc[fillstyle=solid,fillcolor=lightgray](45,-15){39.2}{23}{157}
      \psarc[fillstyle=solid,fillcolor=white](46,-15){38.3}{24}{156}
      \psline[doubleline=true](0,0)(100,0)
      \psline[doubleline=true,doublecolor=black](10,0)(25,0)
      \psline[doubleline=true,doublecolor=black](40,0)(60,0)
      \pscurve[linewidth=.5mm](15,0)(15,7)(30,4)(27,12)(50,9)(45,4)(50,0)
      \put(16,-5){$a$}
      \put(31,-5){$b$}
      \put(48,-5){$c$}
      \put(69,-5){$d$}
      \put(3,7){$e_0$}
      \psline{->}(5,6)(9,2)
      \put(86,6){$v_0$}
      \psline{->}(85.5,5.5)(81.5,1.5)
    \endpspicture
    \caption{\label{fig:boltzmann} Conditioned on the event $B$ (existence
      of the triangle containing $e_0$ and $v_0$), the finite region
      contains a sample of $\tP$. A crossing in the Boltzmann map (shown)
      is expressed as an event in the half plane.}
  \end{center}
  \end{figure}

  Define also the event $A$, to be the event that in the half planar map
  there is a black crossing. Then
  \[ 
    \tP(\tA) = \P(A|B) = \frac{\P(A)\P(B|A)}{\P(B)},
  \]
  and we proceed to evaluate these probabilities. $\P(B)$ is the simplest:
  the event $B$ involves a single step of the exploration process, and
  $\P(B) = p_{m-1}$. The event $A$ has probability $Q_{a,b,c}$, which is
  expressed in \thmref{3segs} in terms of a pair of random walks on 
  $\Z$.

  It remains to find $\P(B|A)$. To this end we re-examine the exploration
  process that was used to determine occurrence of $A$. Suppose some part
  of the map has been revealed, and it is such that both $e_0$ and $v_0$
  are still on the boundary of the infinite unknown component. The distance
  between them along the boundary of the unknown component is no longer
  $m-1$. If the new distance is $l$, then the probability of $B$
  conditioned on the known part of the map is $p_l$.

  Suppose the map is explored along the percolation interface terminating
  at the edge between the segments of lengths $a$ and $b$. For $A$ to
  occur, the interface must hit the second black segment without hitting
  either of the infinite white segments first. Note that as long as the
  interface does not hit the left infinite white segment, $e_0$ is on the
  boundary of the infinite component. Similarly, $v_0$ is in the boundary
  as long as the interface does not touch any vertex to its right. In
  particular, if $A$ occurs, then at the time a black crossing is found
  both $e_0$ and $v_0$ are still in the boundary of the infinite unknown
  component.

  Recall that the lengths at time $t$ of the black and white segments on
  either side of the explored interface along the updated boundary are
  given by random walks $S,S'$. As long as these are positive, the distance
  along the current boundary between $e_0$ and $v_0$ is given by
  $S_t+S'_t+c+d-1$, as there are $S_t$ black vertices, $S'_t$ white
  vertices, $c$ black vertices, and finally $d$ white ones, the last of
  which is $v_0$. For $A$ to occur it is necessary that $\tau = T'_-<T_-$
  and $|S'_\tau| < c$. Conditioned on this, the distance from $e_0$ to
  $v_0$ is now given by $l=S_\tau + S'_\tau + c+d-1$.

  It follows that
  \[
    \P(B|A) = \E_{a,b} \left( p_l \, \big| \, 
                          \big\{ \tau = T'_- \big\} \cap
                          \big\{ |S'_\tau| < c \big\} \right).
  \]
  Combining the three probabilities concludes the proof.
\end{proof}

\subsection{Simple exploration}

It is instructive to consider direct exploration of a Boltzmann
triangulation as a means to evaluate crossing probabilities. Let us
consider the following scenario: the boundary of an $m$-gon is partitioned
into three segments. The first, of length $a$, is black. That is followed
by $b$ white vertices, and finally there are $c$ uncolored vertices. The
interior is filled with a Boltzmann triangulation with critical percolation
on internal vertices. From topological considerations, there is a unique
uncolored boundary vertex that is adjacent to both the black and the white
clusters of the colored boundary segments. Denote this vertex by $w$. We
are interested in the distribution of this vertex along the uncolored
segment.

One can identify $w$ by following the interface between the black and
white clusters across the disc, starting at the edge connecting the two
colored segments. The interface must eventually reach a triangle that
contains one of the uncolored boundary vertices. The first such vertex
encountered is necessarily $w$ (at which point the interface cannot be
followed any longer without coloring the vertex).

As the exploration proceeds, at each step some part of the map has been
revealed. The unknown part consists of a number of discs containing free
triangulations, of which all but one have monochromatic boundaries. The
final disc's boundary contains the uncolored segment, as well as parts of
the white and black percolation clusters of the original boundary segments.
Once the uncolored segment is reached and split into more than one boundary
component, $w$ is determined.

The distribution of $w$ along the uncolored segment is determined by the
lengths of the three segments, so it is natural to consider the evolution
of the segment lengths as the map is explored. Let $A_n,B_n$ be the lengths
of the black and white segments after $n$ exploration steps, assuming that
the exploration has not reached the uncolored segment by that time. We see
that $(A_n,B_n)$ is a Markov chain on $(\Z^+)^2$. If $M = M_n = A_n+B_n+c$.
denotes the length of the boundary after $n$ steps, then the transition
probabilities are:
\[
  (A_{n+1},B_{n+1}) = \begin{cases}
    \begin{array}{c} (A_n+1,B_n\phantom{ +1})\\ 
                     (A_n\phantom{ +1 },B_n+1)\\ \end{array} &  
      \text{with prob. $\frac{Z_{M+1}}{2\alpha Z_M}$ each,} \\
    \begin{array}{c} (A_n-k,B_n\phantom{-k})\\ 
                     (A_n\phantom{-k},B_n-k)\\ \end{array} &  
      \text{with prob. $\frac{Z_{M-k}Z_{k+1}}{Z_M}$ each.} \\
  \end{cases}
\]
Of course, the last two cases hold only for $k<A_n$ and $k<B_n$
respectively. Otherwise --- accounting for the remaining probability ---
the exploration hits the uncolored segment at $w$. The probability of
this happening at step $n+1$ and that $w$ is the uncolored vertex at
distance $k$ from the black segment is
\[
  \frac{Z_{A_n+k+1} Z_{B_n+c-k}}{Z_M}.
\]

\medskip

A direct attempt to scale this process gives only some information about
the properties of its scaling limit, but not a proof that the limit exists.
The natural scaling is to start the process with values $A_0 = \lambda a,
B_0 = \lambda b$, and have $\lambda c$ uncolored vertices (all rounded to
integers), then set
\[
  X^{(\lambda)}_t = \lambda ^{-1} A_{\lambda^{3/2} t}, \qquad
  Y^{(\lambda)}_t = \lambda ^{-1} B_{\lambda^{3/2} t}.
\]
The natural hope is that $(X^{(\lambda)}_t,Y^{(\lambda)}_t)$ converge to
some process $(X_t,Y_t)$ as $\lambda\to\infty$, but we have no direct proof
of this. An indirect proof based on the half plane exploration process
appears in \cite{cardy2}. However, some properties of the scaled process
can be observed from the discrete process.

\begin{thm}
  Assume that the scaling limit of $(X^{(\lambda)}_t,Y^{(\lambda)}_t)$ with
  $\lambda c$ uncolored vertices exists. Then the limit $(X_t,Y_t)$ is
  Markovian, makes jumps from $(X,Y)$ to $(X'+dX,Y)$ for any $0<X'<X$ at
  rate
  \[
    \gamma \left(\frac{(X-X')(X'+Y+c)}{X+Y+c}\right)^{-5/2} dX,
  \]
  for some constant $\gamma$. A similar formula holds for jumps to
  $(X,Y'+dY)$. The limit process is stopped by hitting an uncolored vertex
  at distance $z+dz$ from the black segment, for $0<z<c$, at rate
  \[
    \gamma \left(\frac{(X+z)(Y+c-z)}{X+Y+c}\right)^{-5/2} dz,
  \]
\end{thm}

\begin{proof}
  The discrete process makes a jump from $(A,B)$ to $(A-k,B)$ at a given
  step with probability $Z_{M-k}Z_{k+1}/Z_M$ where $M=A+B+c$. If the
  process is started at $(\lambda a,\lambda b)$ with $\lambda c$ uncolored
  vertices, the probability of a jump of size $\lambda k$ at a given step is
  \[
    \frac{Z_{\lambda(M-k)}Z_{\lambda k+1}}{Z_{\lambda M}}
    \approx 9\gamma' \lambda^{-5/2} \left(\frac{k(M-k)}{M}\right)^{-5/2},
  \]
  where we used the asymptotic formula $Z_n \approx \gamma' 9^n n^{-5/2}$.
  Since there are $\lambda dX$ vertices in the interval $[X',X'+dX]$ and
  since time is scaled by a further factor of $\lambda^{3/2}$, substituting
  $M=X+Y+c$ and $k=X-X'$ yields the first formula. A similar formula for
  jump rates of $Y$ follows from symmetry of the process.

  To estimate termination rates, apply the same asymptotics to the
  probability of hitting a particular vertex in the uncolored segment.
\end{proof}


\begin{thebibliography}{10}

\bibitem{ADJ}
J.~Ambj{\o}rn, B.~Durhuus, and T.~Jonsson.
\newblock {\em Quantum Gravity, a Statitstical Field Theory Approach}.
\newblock Cambridge Monographs on Mathematical Physics, 1997.

\bibitem{AmWa}
J.~Ambj{\o}rn and Y.~Watabiki.
\newblock Scaling in quantum gravity.
\newblock {\em Nucl. Phys. B}, 445(1):129--142, 1995.

\bibitem{cardy2}
O.~Angel.
\newblock Scaling of percolation on infinite planar maps, {II}.
\newblock in preparation.

\bibitem{UIPT2}
O.~Angel.
\newblock Growth and percolation on the uniform infinite planar triangulation.
\newblock {\em Geom. and Func. Anal.}, 13(5):935--974, 2003.
\newblock arXiv:math.PR/0207153.

\bibitem{UIPT1}
O.~Angel and O.~Schramm.
\newblock Uniform infinite planar triangulations.
\newblock {\em Comm. in Math. Phys.}, 241:191--213, 2003.
\newblock arXiv:math.PR/0207153.

\bibitem{BFSS}
C.~Banderier, P.~Flajolet, G.~Schaeffer, and M.~Soria.
\newblock Random maps, coalescing saddles, singularity analysis, and airy
  phenomena.
\newblock {\em Rand. Struc. Alg.}, 19(3-4):194--246, 2001.

\bibitem{BRW}
E.~A. Bender, B.~L. Richmond, and N.~C. Wormald.
\newblock Largest 4-connected components of 3-connected planar triangulations.
\newblock {\em Rand. Struc. Alg.}, 7(4):273--285, 1995.

\bibitem{bert}
J.~Bertoin.
\newblock {\em L\'evy Processes}.
\newblock Cambridge University Press, 1998.

\bibitem{ChSc}
P.~Chassaing and G.~Schaeffer.
\newblock Random planar lattices and integrated super-brownian excursion.
\newblock {\em Prob. Th. and Rel. Fields}, 128(2):161--212, 2004.
\newblock arXiv:math.CO/0205226.

\bibitem{dupl2}
B.~Duplantier.
\newblock Random walks and quantum gravity in two dimensions.
\newblock {\em Phys. Rev. Lett.}, 81(25):5489--5492, 1998.

\bibitem{dupl1}
B.~Duplantier.
\newblock Random walks, polymers, percolation, and quantum gravity in two
  dimensions.
\newblock {\em Phys. A}, 263(1--4):452--465, 1999.

\bibitem{KPZ}
V.~G. Knizhnik, A.~M. Polyakov, and A.~B. Zamolodchikov.
\newblock Fractal structure of 2d-quantum gravity.
\newblock {\em Mod. Phys. Lett. A}, 3:819--826, 1998.

\bibitem{MaMo}
J.~F. Marckert and A.~Mokkadem.
\newblock Limit of normalized quadrangulations: the brownian map.
\newblock arXiv:math.PR/0403398.

\bibitem{RiWo}
B.~L. Richmond and N.~C. Wormald.
\newblock Random triangulations of the plane.
\newblock {\em Euro. J. Comb.}, 9(1):61--71, 1988.

\bibitem{Smir}
S.~Smirnov.
\newblock Critical percolation in the plane.
\newblock www.math.kth.se/ {$\sim$}stas/papers/percol.ps.

\bibitem{Tutte1}
W.~T. Tutte.
\newblock A census of planar triangulations.
\newblock {\em Canad. J. Math.}, 14:21--38, 1962.

\bibitem{Tutte2}
W.~T. Tutte.
\newblock A census of planar maps.
\newblock {\em Canad. J. Math.}, 15:249--271, 1963.

\bibitem{Wern}
W.~Werner.
\newblock Random planar curves and schramm-loewner evolutions.
\newblock Lecture notes for the St. Flour summer school, 2002.
\newblock Springer L.N. Math., 1840, arXiv:math.PR/0303354.

\bibitem{zolo}
V.~Zolotarev.
\newblock {\em One Dimensional Stable Distributions}, volume~65 of {\em Transl.
  of Math. Monographs}.
\newblock American Mathematical Society, 1986.

\end{thebibliography}
\end{document}